\documentclass[a4paper,10pt,english]{amsart}
\usepackage{babel}
\usepackage{amsmath,amsfonts,amssymb,amsthm,mathrsfs,enumitem}
\usepackage[pagebackref=true]{hyperref}
\usepackage{cleveref}
\usepackage[a4paper,margin=3.5cm]{geometry}
\usepackage{setspace}
\allowdisplaybreaks
\theoremstyle{plain}
\newtheorem{theorem}{Theorem}[section]
\newtheorem*{theorem*}{Theorem}
\newtheorem{proposition}[theorem]{Proposition}
\newtheorem*{proposition*}{Proposition}
\newtheorem{lemma}[theorem]{Lemma}
\newtheorem{corollary}[theorem]{Corollary}

\newtheorem*{claim}{Claim}
\theoremstyle{definition}
\newtheorem{definition}[theorem]{Definition}

\newtheorem{convention}[theorem]{Convention}

\newtheorem{remark}[theorem]{Remark}

\newtheorem*{question}{Question}

\newtheorem{notation}[theorem]{Notation}
\newcommand{\ha}{\ensuremath{\mathbf{HA}}} 
\newcommand{\qa}{\ensuremath{\mathbf{Q}}} 

\newcommand{\pa}{\ensuremath{\mathbf{PA}}} 
\newcommand{\mpr}{\ensuremath{\mathbf{MP}}}  
\def\vp{\varphi}
\def\con{\mathit{Con}}
\def\RFN{\mathit{Urf}}
\def\rfn{\mathit{Lrf}}
\def\timess{\cdot}
\newcommand{\pr}{\mathit{Pr}}
\newcommand{\prf}{\mathit{Prf}}
\newcommand{\ax}{\mathit{Ax}}
\newcommand{\fm}{\mathit{Fml}}
\newcommand{\st}{\mathit{Snt}}

\newcommand{\seq}{\mathit{Seq}}
\newcommand{\infer}{\mathit{Rule}}
\newcommand{\lh}{\mathit{lh}}
\newcommand{\num}{\mathit{Nm}}
\newcommand{\sub}{\mathit{Sb}}
\newcommand{\subb}{\mathit{S}}
\newcommand{\numm}{\mathit{N}}
\newcommand{\N}{\mathbb{N}}

\newcommand{\ko}{\mathcal{O}}
\newcommand{\W}{\mathcal{P}}
\newcommand{\pair}[1]{\langle #1\rangle}

\makeatletter     % cutoff subtraction
\newcommand{\prc}{\mathbin{\mathpalette\prc@inner\relax}}
\newcommand{\prc@inner}[2]{%
	\vbox{\offinterlineskip\m@th
		\ialign{%
			##\cr
			\hidewidth\raisebox{-1.5\height}[0pt][0pt]{$#1.$}\hidewidth\cr
			$#1-$\cr
		}%
	}%
}
\makeatother
\newcommand{\imp}{\rightarrow}

\newcommand{\Imp}{\Rightarrow}
\newcommand{\biimp}{\leftrightarrow}

\title[]{Iterating reflection over intuitionistic arithmetic}
\author[Frittaion]{Emanuele Frittaion}
\address{School of Mathematics,	University of Leeds, UK}
\email{}
\thanks{}

\begin{document}
	
	
\begin{abstract} 
In this note, we investigate iterations of consistency, local and uniform reflection over $\ha$ (Heyting Arithmetic). In the case of uniform reflection, we give a new proof of Dragalin's extension  of Feferman's completeness theorem to $\ha$, drawing on Rathjen's proof of Feferman's classical result (see \cite{PRS25}).
%Unlike previous works, this paper does not contain  an adaptation of Feferman's proof to the intuitionistic case. Similarly,  Rathjen's proof  does not translate to intuitionistic logic, but a  simple  modification will do the trick. 
\end{abstract}
	
\maketitle

\section{Introduction}
	
Let $T$ be a recursively enumerable theory in the language of first-order arithmetic. Throughout, we work in intuitionistic logic. In particular, classical $\pa$ (Peano Arithmetic) may be viewed as $\ha + \mathbf{LEM}$, where $\mathbf{LEM}$ comprises all instances of the law of excluded middle.  
	
Informally, we wish to iterate reflection along the ordinals by setting
\begin{align*}
	T_{0}&= T,  \\
	T_{\alpha+1}&=  (T_\alpha)', \\
	T_{\lambda} &= \bigcup_{\alpha\, <\, \lambda} T_\alpha, \qquad \lambda \text{ limit},
\end{align*}
where $T’$ is obtained from $T$ by adjoining  either $\con(T)$ (consistency), $\rfn(T)$ (local reflection) or $\RFN(T)$ (uniform reflection). This iteration is possible so long as the ordinals are effectively presented, and in particular countable.  
	
Feferman \cite{F62} considers sequences of recursively enumerable theories $(T_d)_{d \in \mathcal{O}}$ indexed by ordinal notations in Kleene’s $\mathcal{O}$ --- what he calls \emph{transfinite recursive progressions}. We will simply say that a formula $\varphi$ is provable in an iteration over $T$ if it is provable in $T_d$ for some $d \in \mathcal{O}$.

\begin{theorem}[Feferman's completeness theorem; {\cite[Thm.\  5.13]{F62}}]
	% {\cite[Thm.\  5.13 p.\ 308]{F62}}]
All true sentences of arithmetic are provable in iterations of uniform reflection over $\pa$.
\end{theorem}   
%Moreover, the following holds. 
\begin{theorem}[{\cite[Thm.\ 4.1, Thm.\ 4.5]{F62}}]
	% [{\cite[Thm.\ 4.1 p.\ 287, Thm.\ 4.5 p.\ 289]{F62}}]
	The sentences provable in iterations of either consistency or local reflection over $\pa$ coincide with those provable in $\pa$ $+$ all true $\Pi_1$ sentences. %joined with all true $\Pi_1$ sentences.	
\end{theorem}

The idea of iterating consistency --- and the fact that such iterations capture all true $\Pi_1$ sentences --- goes back to Turing~\cite[p.\ 210]{T39}.

In this note we prove the corresponding results for $\ha$ with respect to iterations of reflection along Kleene's $\ko$. 

\begin{theorem}[Dragalin's theorem {\cite[Thm.\ 1]{Dra}}; see Theorem \ref{uniform}]\label{Dragalin}
	The sentences provable in iterations of uniform reflection over $\ha$ are precisely those provable  in  $\ha$  extended by the recursive $\omega$-rule. This holds for every recursively enumerable extension of $\ha$. 
\end{theorem}
The situation for consistency and local reflection mirrors the classical case.
\begin{theorem}[see Theorem \ref{local}]
	The sentences provable in iterations of either consistency or local reflection over $\ha$ coincide with those provable in $\ha$  $+$ all true $\Pi_1$ sentences. % joined with all true $\Pi_1$ sentences. 
	In fact, this holds for every recursively enumerable extension of $\ha$.
\end{theorem}

We will establish these results for $\mathbf{HA}$ itself, as the proofs carry over \emph{mutatis mutandis} to all r.e.\ extensions.
	
\subsection{Discussion and Related Work}
	
That uniform reflection over $\ha$ behaves as it does should come as no great surprise since the true core of Feferman’s completeness theorem is the following.
\begin{theorem}[cf.\ {\cite[Thm.\ 5.10]{F62}}]\label{classical}
	Every sentence provable in $\pa$ with the aid of the recursive $\omega$-rule can be proved in some iteration of uniform reflection over $\pa$. 
	%A sentence is provable in an iteration of uniform reflection over $\pa$ if it has a recursive $\omega$-proof.	
\end{theorem}

In fact, Feferman relies on results of Shoenfield~\cite{KSW60,S59}, who showed that every true sentence of arithmetic is provable in $\mathbf{PA}$ extended by the recursive $\omega$-rule.

Dragalin~\cite{Dra} gives only a sketch of the proof of \Cref{Dragalin}.  Another proof  is outlined in Sundholm \cite[Ch.\ 7]{S83}. Both note that the main obstacle in trying to extend Theorem \ref{classical}   to $\ha$ is that  Feferman's key \cite[Thm.\ 4.9]{F62}  depends crucially on \emph{classical} reasoning --- specifically, the use of $\Sigma_1$-$\mathbf{LEM}$,  excluded middle for $\Sigma_1$ formulas. 
	
Dragalin solves this by observing --- or claiming --- that a suitable weakening of Feferman's \cite[Th.\ 4.9]{F62} suffices, whereas Sundholm notes that Derevyankina's approach \cite{D74} to iterated reflection, based on a more flexible  notation system for constructive ordinals, goes through in intuitionistic logic.

In the present note,  we retain Feferman's framework and offer a new, self-contained proof of Dragalin's theorem (see Section \ref{tech} below for more details).

Dragalin's short paper~\cite{Dra}  contains two additional results  worth mentioning. 
	
\begin{definition}
A formula is \emph{negative} if it has no occurrences of $\lor$ or $\exists$.  
It is \emph{quasi negative} if it consists of a block of quantifiers followed by a negative matrix.
\end{definition}
\begin{theorem}[{\cite[Thm.\ 3]{Dra}}]\label{quasi negative}
If a quasi negative formula $\vp$ is realizable in the sense of Kleene's recursive realizability, then $\vp$ has an intuitionistic recursive $\omega$-proof. 
\end{theorem}
The proof of this theorem is based on a primitive recursive bottom-up search for an $\omega$-proof of $\vp$ within an intuitionistic, cut-free $\omega$-calculus with two infinitary $\omega$-rules: one for introducing  $\exists$ in the antecedent and the other for introducing $\forall$ in the succedent.   Now, for negative $\vp$, a systematic proof search can be arranged so that the resulting tree is an $\omega$-proof if and only if $\vp$ is true; notice that a negative formula is true if and only if it is realizable (cf.\ \cite[Lemma 3.2.11]{T73}). For quasi negative $\vp$, a given realizer allows to effectively deal with the  front quantifiers.    In particular, the following holds. 	
\begin{corollary}[Completeness for the negative fragment]
	All true negative sentences are provable in iterations of uniform reflection over $\ha$. In particular, the G\"odel-Gentzen translation of every true sentence is likewise provable in such an iteration.
\end{corollary}	
This also follows because the double negative translation extends  to both $\mathbf{PA}$ and $\mathbf{HA}$ equipped with the recursive $\omega$-rule, and a negative formula is provably equivalent to its own translation within $\ha$.
	
The next result gives a characterization of realizability in terms of iterated uniform reflection. Let $\mpr$ be  Markov's principle, that is, the scheme
\[   \neg\neg \vp\imp  \vp, \] 
for $\vp\in\Sigma_1$, and  let  $\mathbf{ECT}_0$ be the extended Church's thesis scheme 
\[ \forall x\, (\alpha\imp \exists y\, \vp(x,y))\imp \exists e\, \forall x\, (\alpha\imp \exists y\, (\{e\}(x)\simeq y\land \vp(x,y))), \]
for $\alpha$ almost negative.\footnote{A formula is called almost negative if it built from $\Sigma_1$ formulas by means of $\land, \imp, \forall$.}
\begin{theorem}[{\cite[Thm.\ 4]{Dra}}]\label{real}
	A sentence  is  realizable if and only if  it is provable in $\ha+\mathbf{ECT}_0+\mpr$ augmented by the recursive $\omega$-rule.\footnote{All results in Dragalin's paper are claimed without a proof. For a proof of Theorem \ref{real},  see Sundholm \cite[Ch.\ 7:18]{S83}. }
\end{theorem}

Kurata~\cite{K65} appears to be the first to study iterated reflection over $\mathbf{HA}$.  
While he does not discuss $\omega$-provability, he investigates iterations of uniform reflection over both $\mathbf{HA}$ and $\mathbf{HA}+\mathbf{MP}$.  The following results  from \cite{K65} deserve attention.\footnote{In his review of Kurata's paper, Feferman \cite{F73review} points out a minor  error and suggests a quick fix.} In each case,  $T_d$  denotes the $d$-th iteration of uniform reflection over $T$.
\begin{theorem}[{\cite[Thm.\ 1 p.\ 149, Prop.\ 5.5 p.\ 160, Thm.\ 3 p.\ 161]{K65}}]
	The following hold:
	\begin{enumerate}[label=$(\arabic*)$,leftmargin=7mm]
		\item if $\vp$ is provable in  $\ha_d$  for some $d\in\ko$, then $\vp$ is  realizable; 
		\item  if $\vp$ is provable in   $(\ha+\mpr)_d$ for some $d\in\ko$, then  the G\"odel's Dialectica interpretation  $\vp^D$ of $\vp$ is true in  Kreisel's model $\mathbf{HRO}$ of the hereditarily recursive operations;
		\item  if $\vp$ is in prenex form, then $\vp$ is provable in  $(\ha+\mpr)_d$ for some $d\in\ko$ if and only if  $\vp$ is recursively true.
	\end{enumerate}
\end{theorem}

A prenex sentence $\vp$ is recursively true  if $\mathbf{HRO}\models \vp^S$, where $\vp^S$ is the normal Skolem form of $\vp$;
cf.\ \cite[Ch.\  XV, Sec.\  79, p.\ 465]{K71}. Note that a sentence in prenex form is recursively true if and only if it is realizable (see \cite[Ch.\  XV, Sec.\  82, p.\ 516]{K71}).  It is not difficult to see that any such sentence has an intuitionistic recursive $\omega$-proof. By Dragalin's theorem, one obtains the following interesting corollary.

\begin{theorem}
	Iterations of uniform reflection over $\ha$ and over $\ha+\mpr$ prove the same prenex sentences. 
	%$\bigcup_{d\in\ko} \ha_d$ and $\bigcup_{d\in\ko} (\ha+\mpr)_d$ prove the same prenex sentences.
\end{theorem}

\subsection{Motivation and Open Questions}
	
Our interest  stems from a wish to better understand the intuitionistic recursive $\omega$-rule.  
Since the recursive $\omega$-rule over $\mathbf{PA}$ is complete for true sentences of arithmetic, it is natural to ask how it fares against the unrestricted $\omega$-rule when the underlying logic is intuitionistic.
	
It is well known that
\begin{enumerate}
		\item every true arithmetical sentence has an intuitionistic $\omega$-proof (albeit one that is far from constructive, as it can only be found by using the set of true sentences as an oracle), and
		\item any sentence with an intuitionistic recursive $\omega$-proof is realizable; cf.\ \cite[p.\ 91]{L77}.
\end{enumerate}
	
Since there are true sentences that are not realizable, it follows that the recursive $\omega$-rule is strictly weaker than the full $\omega$-rule in the intuitionistic setting. Notice that each instance of excluded middle for $\Sigma_1$ sentences --- and hence Markov’s principle for closed ones --- is realizable and  has an intuitionistic recursive $\omega$-proof, though neither the realizer nor the proof can be found effectively.  
	
Admittedly, the  characterization in terms of reflection principles offers little deeper insight into the intuitionistic recursive $\omega$-rule than the rule itself, yet it makes for a nice companion to the classical analysis of iterated uniform reflection.
	
\begin{question}
Is every realizable true sentence provable in an iteration of uniform reflection over $\mathbf{HA}$?  
If not, is there a notion of intuitionistic truth that renders such iterations complete?  
And what of $\mathbf{HA}+\mathbf{MP}$?
\end{question}
	
We claim that $\mpr$ is not provable in iterations of uniform reflection over $\mathbf{HA}$.  
This gives a counterexample --- a realizable, true sentence not derivable in that way.  Therefore, the answer to the first question is negative.  Kreisel (unpublished; see \cite{F73review} and \cite[Sec.\ 7:17]{S83}) observed that the closure of $\ha$ under Markov's rule extends to every iteration of uniform reflection  $\ha_d$ for $d\in \ko$, and thus to  $\ha$  equipped with the recursive $\omega$-rule. This suggests that any notion of truth that seeks to achieve the desired completeness with respect to recursive $\omega$-provability over $\ha$ cannot be compositional.   Specifically, the following Tarski's clause for implication
\[  T(\vp)\imp T(\psi) \text{ implies } T(\vp\imp\psi)\]
must be dropped.

Finally, since all our proofs rely squarely on intuitionistic logic (see the implications (2)$\Rightarrow$(3) of Theorem \ref{local} on local reflection and (1) $\Imp$ (2) of Theorem \ref{uniform} on uniform reflection), it is natural to ask:
	
\begin{question}
Do any of these results   hold for minimal logic, that is, intuitionistic logic without the principle of explosion (ex falso quodlibet)?
\end{question}

\subsection{Technical Note}\label{tech}
	
Our proof of Dragalin’s theorem (see \Cref{uniform}) builds upon Rathjen's proposal~\cite{R18,R22} of using Sch\"utte’s canonical search trees \cite{Schu56}  to prove Feferman’s completeness theorem for uniform reflection.  
For any sentence $\varphi$, one can construct its canonical tree $S(\varphi)$ (its \emph{Stammbaum}) in a primitive recursive fashion. The result is that $S(\varphi)$ is a classical $\omega$-proof of $\varphi$ if and only if $\varphi$ is true.   In the intuitionistic case, we cannot replicate a similar scenario; nonetheless --- and here lies our key idea --- we can construct a primitive recursive tree $S(a)$ for any number $a$ such that $S(a)$ is an intuitionistic $\omega$-proof of $\varphi$ (possibly with repetitions) if and only if $a$ codes an intuitionistic recursive $\omega$-proof of $\varphi$. Elaborating this idea yields \Cref{uniform}.   
	
A further insight, due to Rathjen, is the use of L\"ob's theorem to formalize certain arguments within $\mathbf{PA}$. Although L\"ob’s theorem does not hold wholesale for $\mathbf{HA}$, its restriction to  $\Pi_2$ formulas does, and this will allow us to prove --- within $\mathbf{HA}$ itself --- several facts  that {\em prima facie} require induction on Kleene's $\ko$.

\subsection{Outline} In Section \ref{prelim} we gather all necessary ingredients. The point of this lengthy review is to fix terminology and notation, and  ensure that $\ha$ has all  the required apparatus.    Readers may skip ahead to Sections \ref{local sec} and \ref{uniform sec}, where the main proofs are presented, referring back to Section \ref{prelim}  only as needed.

\section{Preliminaries}\label{prelim}

Let $\ha$ be a standard axiomatization of Heyting Arithmetic with only one relation symbol~$=$ for equality  and function symbols for  all primitive recursive functions;  cf.\ Troelstra \cite[Ch.\  1, Sec.\  3]{T73}. Negation $\neg\vp$ is defined as $\vp\imp\bot$.

\subsection{Arithmetization}

If $F\colon\N^k\to\N$ is a primitive recursive function, we denote by $F^\circ$ the corresponding function symbol.\footnote{Indeed, $F^\circ$ corresponds to a given {\em description} of the primitive recursive function $F$.} In particular,
\[ \ha\vdash F^\circ(\bar n_1,\ldots,\bar n_k)=\overline{F(n_1,\ldots,n_k)}. \] 
Moreover, if $R\subseteq\N^k$ is a primitive recursive relation, we denote by $R^\circ(x_1,\ldots,x_k)$ the formula $F^\circ(x_1,\ldots,x_k)=0$, where $F$ is the {\em representing} function of $R$.\footnote{The function that takes value $0$ when $(n_1,\ldots,n_k)\in R$ and value $1$ otherwise; cf.\ Kleene \cite[p.\ 8]{K71}.} In particular, 
\begin{align*}
&	(n_1,\ldots,n_k)\in R  \quad\text{iff} \quad \ha\vdash R^\circ(\bar n_1,\ldots,\bar n_k), \\
&	(n_1,\ldots,n_k)\notin R \quad \text{iff}\quad  \ha\vdash \neg R^\circ(\bar n_1,\ldots,\bar n_k). 
\end{align*}
In other words, primitive recursive functions and relations are representable\footnote{In the sense of Shoenfield; see \cite[Ch.\  6, Sec.\  7]{Sbook}.} in $\ha$.

\begin{convention}
We typically drop the $\circ$ notation when the context makes it clear whether we are talking about the primitive recursive function (relation) or its corresponding symbol. 
\end{convention}

We can define $x<y$ by $S(x)\prc y=0$, where $\prc$ is the primitive recursive cutoff subtraction. Bounded quantifiers are then defined as usual. 
\begin{definition}[{cf.\ \cite{V82}}]
The $\Delta_0$ formulas are defined by the following clauses:
\begin{itemize}[leftmargin=5mm]
	\item atomic formulas (including $\bot$) are $\Delta_0$;
	\item if $\vp$ and $\psi$ are $\Delta_0$, then $\vp\land\psi$, $\vp\lor\psi$ and $\vp\imp\psi$ are $\Delta_0$;
	\item if $\vp$ is $\Delta_0$ and $t$ is a term not containing $x$, then $\forall x<t\, \vp(x)$ and $\exists x<t\, \vp(x)$ are $\Delta_0$.
\end{itemize}
Let
\begin{align*}
\Sigma_1&=\{\exists x\, \vp(x)\mid \vp\in\Delta_0\}, &
\Pi_1&=\{\forall x\, \vp(x)\mid \vp\in\Delta_0\}, &
\Pi_2 &=\{\forall x\, \vp(x)\mid \vp\in \Sigma_1\}.
\end{align*}
The $\Sigma$ formulas are defined by the following clauses:
\begin{itemize}[leftmargin=5mm]
 \item $\Delta_0$ formulas are $\Sigma$;
 \item  if $\vp$ and $\psi$ are $\Sigma$, then $\vp\land\psi$ and $\vp\lor\psi$ are $\Sigma$;
 \item if $\vp$ is $\Delta_0$ and $\psi$ is $\Sigma$, then $\vp\imp \psi$ is $\Sigma$;
 \item if $\vp$ is $\Sigma$ and $t$ is a term not containing $x$, then $\forall x<t\, \vp(x)$ is $\Sigma$;
 \item if $\vp$ is $\Sigma$, then $\exists x\, \vp(x)$ is $\Sigma$.
\end{itemize}
The $\Pi$ formulas are defined dually by switching $\exists$ and $\forall$ in last two clauses. 
\end{definition}

By a routine induction on the build-up of a formula one can prove the following.

\begin{proposition}[normal form]
For every $\Delta_0$ formula $\vp(x_1,\ldots,x_k)$ there is a primitive recursive function $F$ such that 
\[ \ha\vdash \vp(x_1,\ldots,x_k) \biimp F(x_1\,\ldots,x_k)=0. \]
For every $\Sigma$ formula $\vp(x_1,\ldots,x_k)$ there is a primitive recursive function $F$ such that 
\[ \ha\vdash \vp(x_1,\ldots,x_k) \biimp \exists x\, F(x_1\,\ldots,x_k,x)=0. \]
For every $\Pi$ formula $\vp(x_1,\ldots,x_k)$ there is a primitive recursive function $F$ such that 
\[ \ha\vdash \vp(x_1,\ldots,x_k) \biimp \forall x\, F(x_1\,\ldots,x_k,x)=0. \]
\end{proposition}
In particular, $\Delta_0$ formulas are decidable in $\ha$, namely,  $\ha\vdash \vp\lor\neg\vp$,  for every $\Delta_0$ formula $\vp$.

\begin{definition}
A formula $\vp$ is  $\Sigma$ in $\ha$ if it is equivalent in $\ha$ to a $\Sigma$ formula. Similarly for other classes of formulas. 
\end{definition}

We  assume a primitive recursive coding of finite sequences. As in Feferman \cite{F60,F62}, we make  no distinction between  expressions  (terms or formulas) and their G\"{o}del numbers. The logical symbols $\land, \lor,\imp, \forall, \exists$ can be treated as (primitive recursive) operations on the natural numbers. We fix a Hilbert-style system for intuitionistic first-order logic with equality. We may safely assume a primitive recursive presentation of the logical axioms and rules.

\begin{notation}
We make use of the following  primitive recursive functions and relations.  The predicate $\seq(x)$ indicates  that $x$ is sequence, $\lh(x)$ is  the length of $x$, $x_i$ is the $i$-th element of $x$ and $x\restriction  i$ is the initial segment of $x$ of length $i$. Let $\st(x)$ and $\fm(x)$ express that $x$ is a sentence and $x$ is a formula, respectively.  The predicate $\ax(v)$ means that $x$ is a logical axiom and $\infer(x,y)$ means that the formula $x$ can be obtained from  formulas in the sequence $y$ by means of a logical rule.  
\end{notation}

\begin{definition}[provability predicate]
Let $\alpha(v)$ be a formula. The formula $\pr_{\alpha(v)}(x)$ is defined as $\exists y\, \prf(x,y)$, where
\[ \prf(x,y) \ =_\text{def} \  \seq(y)\land x=y_{\lh(y)-1}\land \forall i<\lh(y)\, \big(\ax(y_i)\lor \alpha(y_i)\lor \infer(y_i,y\restriction  i)\big).  \]
The consistency formula $\con_{\alpha(v)}$ is defined by $\neg\pr_{\alpha(v)}(\bar\bot)$.
\end{definition}

\begin{remark}
If $\alpha(v)$ is $\Sigma$ (in $\ha$), then $\pr_{\alpha(v)}(x)$ is $\Sigma$ (in $\ha$) and  $\con_{\alpha(v)}$ is $\Pi$ in $\ha$.  Note that the formula $\alpha(v)$ may contain further free variables in addition to $v$.
\end{remark}

\begin{remark}\label{remark}
Let $\alpha(z,v)$ be any formula with displayed free variables, and let $\pr(z,x)$ be the corresponding provability predicate.   Then 
	\[ \ha\vdash \forall v\, (\alpha(z_0,v)\imp \alpha(z_1,v))\imp \pr(z_0,x)\imp \pr(z_1,x).\]
\end{remark}

\begin{definition}[syntactic operations]
There are primitive recursive functions $\num$, $\sub$, $\pr$ and $\con$ such that:
\begin{itemize}[leftmargin=5mm]
		\item $\num(n)$ is the $n$-th numeral $\bar n$ for every number $n$;
	    \item $\sub(x,y,z)$ is the  substitution of $z$  for $y$ in $x$. In general, $\sub(x,y_1,\ldots,y_k,z_1,\ldots,z_k)$ is the simultaneous substitution of  $z_i$ for $y_i$ in $x$;  
		\item $\pr_\alpha$ and $\con_\alpha$ are the provability predicate and the   consistency formula for $\alpha$, respectively, whenever $\alpha(v)$ is a formula and $v$ is a variable. 
\end{itemize}
\end{definition}

%containing solely the free variables  $z$ and $v$.
%1.  \alpha(z, v)  is a formula with z and v as its only free variables.
%2. \alpha(z, v) is a formula having only \( z \) and \( v \) as its free variables.
%3.   \alpha(z, v)  is a formula in which \( z \) and \( v \) are the only variables free.
%4.  \alpha(z, v) is a formula containing solely the free variables \( z \) and \( v \).
%5.  \alpha(z,v) is a formula with only the free variables z and v

Strictly speaking, $\pr_{\alpha}$ is a function of both  $\alpha$ and the variable $v$, not just $\alpha$.   We rely on context and omit $v$. We can  assume that this operation is always well defined and returns a formula whose free variables are those of $\alpha$ except for $v$ and a fresh new variable $x$. The same goes for $\con_\alpha$.  We will not bother with such details in the future, but the reader should be aware of them.

\begin{definition}[syntactic operations  one level down]
There are  primitive recursive functions $\numm$ and   $\subb$ such that: 
\begin{itemize}[leftmargin=5mm]
	\item $\numm(x)$ is the term $\num^\circ(x)$, whenever $x$ is a term;
	\item  $\subb(x,y,z)$ is  the term $\sub^\circ(x,y,z)$, whenever  $x,y,z$ are terms.
\end{itemize} 
Let $\numm^\circ$ and $\subb^\circ$ be the corresponding function symbols. We  use $\dot x$ to denote $\num^\circ(x)$ and $\ddot{x}$ to denote $\numm^\circ(x)$. We  write $\bar\vp(\dot x)$ for $\sub^\circ(\bar\vp,\bar x,\dot x)$, and   $\dot \vp (\ddot{x})$ for  $\subb^{\circ}(\dot \vp,\dot x,\ddot{x})$.
\end{definition}
\begin{remark}
If $\vp(x)$ is a formula containing solely the free variable  $x$, then 
\[  \ha\vdash \dot{\bar \vp}(\ddot{\bar x}) = \overline{\bar\vp(\dot x)}.\] 
\end{remark}

\begin{definition}[defining reflection]
Let $T$ be a theory defined by a formula  $\alpha(v)$ with a single free variable $v$.
Let  
\[\pr(x)\ =_\text{def}\ \pr_{\alpha}(x)  \quad\text{ and } \quad \con(T) \ =_\text{def}\ \neg \pr(\bar\bot).\]
Write $\pr^\circ$ and $\con^\circ$ for $\pr^\circ_{\bar\alpha}$ and $\con^\circ_{\bar\alpha}$, respectively. 
\begin{itemize}[leftmargin=5mm]
\item The scheme  $\rfn(T)$ consists of all sentences of the form 
\[ \pr(\bar \vp)\imp \vp.  \]
\item The scheme	$\RFN(T)$ consists of all sentences of the form 
\[ \forall x\, (\pr(\bar \vp(\dot x)) \imp  \vp(x)). \]
\end{itemize}
\vspace{1mm}
\begin{enumerate}[label=(\alph*), labelindent=0mm,leftmargin=*]
\item The theory $T+\con(T)$ is  defined by $\alpha(v)\lor  v=\con^\circ$. \vspace{2mm}
\item The theory $T+\rfn(T)$ is defined by $\alpha(v)\lor \vartheta(v)$, where
\[ \vartheta(v) \ =_\text{def}\ \exists \vp\, (\st(\vp)\land v=\pr^\circ(\dot \vp)\imp^\circ \vp).\]
\item The theory $T+\RFN(T)$ is defined by $\alpha(v)\lor \vartheta(v)$, where
\[ \vartheta(v) \ =_\text{def} \ \exists \vp\, \exists x\, (\st(\forall^\circ x\, \vp)   \land v=\forall^\circ x\, (\pr^\circ(\dot \vp(\ddot x))\imp^\circ  \vp(x))). \]
\end{enumerate}
\end{definition}

\begin{remark}
If $\alpha(v)$ is $\Sigma$, then $T'$ is also $\Sigma$ definable. 
\end{remark}

\subsection{Self-reference}

\begin{theorem}[diagonal lemma]\label{diagonal}
	Let $\vp(u,x_1,\ldots,x_n)$ be any formula. Then there is a formula $\delta(x_1,\ldots,x_n)$ such that 
	\[ \ha\vdash \delta(x_1,\ldots,x_n)\biimp \vp(\bar\delta, x_1,\ldots,x_n).\]
\end{theorem}
\begin{proof}
	For ease of notation, let $n=1$.
	Take the primitive recursive function $D\colon \N\to \N$ such that
	\[ D(\psi(u,x))=\psi(\overline{\psi(u,x)},x). \]
	Let $g$ be the (G\"{o}del number of the) formula  $\vp(D^{\circ}(u),x)$ and let
	$\delta(x)\ =_\text{def}\ \vp(D^{\circ}(\bar g),x)$. 
\end{proof}
\begin{remark}
If $\vp$ is $\Sigma$ (in $\ha$), so is $\delta$.
\end{remark}

Let $T(e,x,z)$ be  Kleene $T$-predicate with output function $U(z)$. Write $\{e\}(x)\simeq y$ for $\exists z\, (T(e,x,z)\land U(z)=y)$ and $\{e\}_z(x)\simeq y$ for $T(e,x,z)\land U(z)=y$. Further abbreviations are $\{e\}(x)\downarrow$ for $\exists y\, (\{e\}(x)\simeq y)$ and $\{e\}_z(x)\downarrow$ for  $T(e,x,z)$. 

\begin{theorem}[s-m-n and recursion theorem; {cf.\ Troelstra \cite[Ch.\ 1, Sec.\ 3.10]{T73} } ]
Let $m,n\geq 1$. There are primitive recursive functions $S^m_n$ and $F_n$ such that
\begin{align*}
	\ha& \vdash \{S^m_n(e,z_1,\ldots,z_m)\}(x_1,\ldots,x_n)\simeq \{e\}(z_1,\ldots,z_m,x_1,\ldots,x_n) \\
\ha& \vdash   \{F_n(e)\}(x_1,\ldots,x_n)\simeq \{e\}(F_n(e),x_1,\ldots,x_n) 
\end{align*}
\begin{remark}
The s-m-n theorem is all we need to prove the recursion theorem.
\end{remark}
\begin{remark}
We use the recursion theorem  mainly to compute indices $e$ of total recursive functions satisfying
\[ \ha \vdash \{\bar e\}(x_1,\ldots,x_n)\simeq G(\bar e,x_1\ldots,x_n),\]
where $G$ is some primitive recursive function. 
Moreover, we repeatedly use  the $S^m_n$ functions to primitive recursively construct indices from indices.  
\end{remark} 

\end{theorem} 

Recall that   Robinson's $\qa$ is a finite theory in the language $0,S,+,\timess$.  In keeping with our  formulation of $\ha$, we extend $\qa$ to include all primitive recursive functions.    We may safely assume that  \[ \ha \vdash \pr_{\qa}(x)\imp \pr_\ha(x), \]
where $\pr_T(x)$  denotes $\pr_{\alpha}(x)$ for a standard primitive recursive presentation $\alpha(v)$ of $T$.

\begin{theorem}[$\Sigma$ completeness; {cf.\ \cite[Thm.\  2.9]{V82}}]\label{completeness}
For every $\Sigma$ formula $\vp(x_1,\ldots,x_k)$:
\begin{itemize}[leftmargin=5mm]
	\item if $\vp(\bar n_1,\ldots,\bar n_k)$ is true, then $\qa\vdash \vp(\bar n_1,\ldots,\bar n_k)$;
    \item  $\ha \vdash \vp(x_1,\ldots,x_k)\imp \pr_{\qa}(\bar \vp(\dot x_1,\ldots,\dot x_k))$.
\end{itemize}
In particular, if $\vp(x_1,\ldots,x_k)$ is $\Sigma$ in $\ha$, then 
\begin{itemize}[leftmargin=5mm]
	\item $\ha\vdash \vp(x_1,\ldots,x_k)\imp \pr_{\ha}(\bar \vp(\dot x_1,\ldots,\dot x_k))$.
\end{itemize}
\end{theorem}

The next theorem  (its restriction to $\Sigma$ formulas) is needed to prove the existence of {\em transfinite recursive progressions} in the sense of Theorem \ref{existence} (see however Remark \ref{weak existence}).
  
\begin{theorem}[soundness of $\qa$; {cf.\ Troelstra \cite[Ch.\ 1, Sec.\ 5.9]{T73}}]\label{soundness}
For every formula  $\vp(x_1,\ldots,x_k)$,
\[ \ha \vdash \pr_{\qa}(\bar\vp(\dot x_1,\ldots,\dot x_k))\imp \vp(x_1,\ldots,x_k).\]
\end{theorem}

%%%%%%%%%%%%%%%%%%%%%%%%%%%%
\begin{notation}
For the rest of the paper, $\alpha(z,v)$ will be a formula 	having only $z,v$ as free variables.   We may use $\pr(z,x)$ and $\con(z)$ to denote $\pr_{\alpha}(x)$ and  $\neg\pr(z,\bar\bot)$, respectively. Similarly, we write $\pr^\circ$ and $\con^\circ$ for $\pr^\circ_{\bar\alpha}$ and $\con^\circ_{\bar\alpha}$, respectively. 
\end{notation}

\begin{theorem}[I am not provable; {cf.\ \cite[Thm.\  2.11 p.\ 272]{F62}}]\label{nu}
Let $\alpha(z,v)$ be a $\Sigma$ formula in $\ha$ having only $z,v$ as free variables such that 
\[ \tag{$\ast$} \ha\vdash \pr_\ha(x)\imp \pr(z,x). \]  
Then there is a formula $\nu(z)$ such that 
\[ \ha \vdash \con(z)\biimp \nu(z) \biimp \neg\pr(z,\bar \nu(\dot z)). \]
In particular, $\nu(z)$ is $\Pi_1$ in $\ha$. 
\end{theorem}

\begin{proof}
By the diagonal lemma, there is a fixed point  $\nu(z)$ of the formula  $\vp(u,z)\ =_\text{def}\ \neg\pr(z,\sub^\circ(u,\dot z))$. Then  
$\ha\vdash \nu(z)\biimp \neg\pr(z,\bar \nu(\dot z))$. In particular, $\nu(z)$ is $\Pi_1$ in $\ha$.   Clearly, $\ha\vdash \nu(z)\imp \con(z)$. The converse direction $\ha\vdash \con(z)\imp \nu(z)$ relies on the fact that \[ \tag{0} \ha\vdash \nu(z)\biimp \forall y\, \psi(z,y), \] for some $\Delta_0$ formula $\psi$.\footnote{In the case of $\pa$ one can prove the contrapositive  $\neg\nu(z)\imp \neg\con(z)$ and freely infer  $\pr(z,\bar\nu(\dot z))$ from $\neg\nu(z)$.} Therefore:
\begin{align*}
\tag{1} \ha& \vdash \neg\psi(z,y)\imp \pr_\ha(\overline{\neg \psi}(\dot z,\dot y))	& \text{ $\Sigma$ completeness (indeed $\Delta_0$)}, \\
\tag{2}  \ha& \vdash \neg\psi(z,y) \imp \pr_\ha(\overline{\neg \forall y\, \psi}(\dot z)) & \text{ (1),}\\
\tag{3} \ha&\vdash \neg\psi(z,y)\imp \pr(z, \overline{\neg \forall  y\, \psi}(\dot z))  &  \text{ (2) and ($\ast$),} \\
\tag{4} \ha& \vdash \pr(z, \overline{\forall z\, (\nu(z) \biimp \forall y\, \psi(z,y)}) & \text{ $\Sigma$ completeness from (0)  and ($\ast$),}\\
\tag{5} \ha& \vdash \neg\psi(z,y)\imp \pr(z, \neg^\circ\, \bar\nu(\dot z ))   &    \text{ (3) and (4),} \\
\tag{6} \ha& \vdash \con(z)\imp \neg\psi(z,y)\imp \neg \pr(z, \overline \nu(\dot z))      & \text{ (5),}\\
\tag{7} \ha&\vdash \con(z)\imp \neg\psi(z,y)\imp \nu(z)   & \text{ (6) and the fixed point property,}\\
\tag{8} \ha&\vdash \neg\psi(z,y)\imp  \neg\nu(z)      &  \text{ pure logic from (0),} \\
\tag{9} \ha&\vdash \con(z)\imp \forall y\, \neg\neg \psi(z,y)        & \text{ (7) and (8) (in minimal logic),}\\
\tag{10} \ha&\vdash \con(z)\imp \forall y\, \psi(z,y)   & \text{ (9) and decidability of $\Delta_0$ formulas.} 
\end{align*}
Hence $\ha\vdash \con(z)\imp \nu(z)$, as desired.
\end{proof}

\begin{theorem}[second incompleteness]\label{second}
Let $\alpha(z,v)$ be as in Theorem \ref{nu}. Then 
\[ \ha\vdash \con(z)\imp \neg\pr(z, \con^\circ(\dot z)). \]
\end{theorem}
\begin{proof}
\begin{align*}
\tag{1} \ha &\vdash \con(z)\imp \neg \pr(z,\bar \nu(\dot z)) & \text{ Theorem \ref{nu},} \\
\tag{2} \ha& \vdash \con(z)\imp \nu(z) & \text{ Theorem \ref{nu},} \\
\tag{3} \ha & \vdash \pr_\ha(\overline{\forall z\, (\con(z)\imp \nu(z))}) & \text{ $\Sigma$ completeness from (2),}\\
\tag{4} \ha&\vdash \pr_\ha(\con^\circ(\dot z)\imp^\circ \bar\nu(\dot z)) & \text{ (3),} \\
\tag{5} \ha&\vdash  \pr(z,\con^\circ(\dot z)\imp^\circ \bar\nu(\dot z)) & \text{ (4) and ($\ast$),} \\
\tag{6} \ha& \vdash \pr(z,\con^\circ(\dot z))\imp \pr(z,\bar\nu(\dot z)) & \text{ (5),} \\
\tag{7} \ha&\vdash \con(z)\imp \neg\pr(z, \con^\circ(\dot z)) & \text{ (1) and (6).} 
\end{align*}
\end{proof}

\begin{theorem}[L\"{o}b's theorem for $\Pi_2$ formulas]\label{lob}
Let $\vp$ be $\Pi_2$ in $\ha$. Then 
\begin{align*}  &\ha\vdash \pr_\ha(\bar\vp)\imp \vp \quad\text{iff} \quad \ha\vdash \vp \end{align*}
\end{theorem}
\begin{proof}
Notice that this holds for any formula in the case of $\pa$. On the other hand, $\pa$ is $\Pi_2$-conservative  over $\ha$ and  the proof of this fact via  negative translation coupled with  Friedman A-translation   formalizes in $\ha$.  The theorem thus follows from L\"{o}b's theorem for $\pa$. 
\end{proof}

\subsection{Kleene's $\ko$}
Let us recall Kleene's  notation system $(\ko,<_\ko)$ for constructive  ordinals. The set  $\ko$ consists of natural numbers which are notations for all recursive ordinals. Moreover, $<_\ko$ is a well-founded partial order on $\ko$ such that  $a<_\ko b$ implies $|a|<|b|$, where $|a|$ is the recursive ordinal denoted by $a$. 

\begin{definition}
The sets $\ko$ and $<_\ko$ are inductively defined by:
\begin{itemize}[leftmargin=5mm]
	\item $0\in\ko$;
	\item if $d\in\ko$, then $2^d\in\ko$ and $d<_\ko 2^d$;
	\item if $e$ is the index of a total recursive function and $\{e\}(n)<_\ko \{e\}(n+1)$ for every $n$, then $3\timess 5^e\in \ko$ and $\{e\}(n)<_\ko 3\timess 5^e$ for every $n$;
	\item if $a<_\ko b<_\ko c$ then $a<_\ko c$.
\end{itemize}
We write $a\leq_\ko b$ for $a<_\ko b \lor a=b$.
\end{definition}

\begin{proposition}[properties of  $<_\ko$]\label{O}
The following properties hold:
\begin{itemize}[leftmargin=5mm]
	\item if $a<_\ko b$ then $a,b\in\ko$;
	\item if $a\in\ko$ then $0\leq_\ko a$;
	\item if $a<_\ko b$ then $b\notin\{0,a\}$;
	\item if $a<_\ko 2^b$ then $a\leq_\ko b$;
	\item if $a<_\ko b$ then $2^a\leq_\ko b$;
	\item  if $a<_\ko 3\timess 5^e$ then $a<_\ko \{e\}(n)$ for some $n$.
\end{itemize}
\end{proposition}

For technical reasons (see Lemma \ref{right}), it will be convenient to define addition as follows.\footnote{The usual clause at limits reads $\{H(a,e)\}(n)\simeq a+_\ko \{e\}(n)$.}
\begin{proposition}
There is a total recursive function $+_\ko$ such that:
\[ a+_\ko b\simeq 
\begin{cases}
	a & b=0 \\
	2^{a+_\ko c} &  b=2^c\\
	3\timess 5^{H(a,e)} &  b=3\timess 5^e,   \text{\ where } \{H(a,e)\}(n)\simeq a+_\ko 2^{\{e\}(n)} \\
	7 & \text{otherwise}
\end{cases}  \] 
\end{proposition}
\begin{proof}
By the recursion theorem.
\end{proof}
\begin{remark}
The function $F$ such that 
$F(a,b)=a+_\ko b$ is primitive recursive. 
\end{remark}
\begin{proposition}[properties of $+_\ko$]\label{plus}
The following properties hold:
\begin{itemize}[leftmargin=5mm]
\item $a,b\in\ko$ iff $a+_\ko b\in\ko$;
\item if $a,b\in\ko$ and $b\neq 0$ then $a<_\ko a+_\ko b$;
\item  $a\in\ko$ and $b<_\ko c$ iff $a+_\ko b<_\ko a+_\ko c$; 
\end{itemize}
\end{proposition}
\begin{definition}
Let $0_\ko=0$ and $(n+1)_\ko= 2^{n_\ko}$. Note that $a+_\ko n_\ko<_\ko a+_\ko m_\ko$ for every $a\in\ko$ and for all $n<m$.
\end{definition}
We use association to the left and thus write $a+_\ko b+_\ko c$ for $(a+_\ko b)+_\ko c$.

\subsection{Transfinite Recursive Progressions}

Let us fix a standard primitive recursive presentation $\alpha_0(v)$ of  $\ha$. Let $\pr_\ha(x)$ be the corresponding provability predicate. 

\begin{definition}[progression and succession formulas]
Let $\alpha(z,v)$ and $\rho(u,z,v)$ be $\Sigma$ formulas. We say that  $\alpha(z,v)$ is a  progression formula based on $\rho(u,z,v)$  if 
	\begin{align*} 
	\ha \vdash&  \ \alpha_0(v)\imp \alpha(z,v), \\	
	\ha \vdash & \  \alpha(0,v) \biimp \alpha_0(v), \\
	\ha \vdash & \  \alpha(2^d,v) \biimp \alpha(d,v)\lor \rho(\bar\alpha,d,v),   \\
	\ha \vdash & \ \alpha(3\timess 5^e,v) \biimp \alpha_0(v)\lor \exists n\, \exists d\,  \big(\{e\}(n)\simeq d\land  \alpha(d,v))\big).
\end{align*}
In this context,  $\rho(u,z,v)$ is called a succession formula.
\end{definition}

\begin{theorem}[existence]\label{existence}
Let $\rho(u,z,v)$ be $\Sigma$. Then there is a $\Sigma$ progression formula $\alpha(z,v)$ based on $\rho(u,z,v)$.  	
\end{theorem}
\begin{proof}
By the diagonal lemma (Theorem \ref{diagonal}), there is a fixed point $\alpha(z,v)$ of the formula $\vp(u,z,v)$ so that $\ha\vdash \alpha(z,v)\biimp \vp(\bar\alpha,z,v)$, where 
\begin{multline*}
\vp(u,z,v)\ =_\text{def}\ (\vp_0(z)\land \alpha_0(v))	\lor \exists d\, (z=2^d\land \vp_1(u,d,v))\lor \exists e\, (z=3\timess 5^e\land \vp_2(u,e,v)),
\end{multline*}
and 
\begin{align*}
\vp_0(z)\ &=_\text{def} \  z=0\lor \forall w<z\, (z\neq 2^w\land z\neq 3\timess 5^w), \\
\vp_1(u,d,v)\ &=_\text{def} \  \pi(u,d,v)\lor \rho(u,d,v),\\
\vp_2(u,e,v)\ &=_\text{def} \  \alpha_0(v)\lor \exists n\, \exists d\, (\{e\}(n)\simeq d\land \pi(u,d,v)),
\end{align*} 
with $\pi(u,d,v)\ =_\text{def}\ \pr_{\qa}(u(\dot d,\dot v))$. Since $\vp(u,z,v)$ is $\Sigma$, so is  $\alpha(z,v)$.  The first item in the definition of progression formula is proved by ordinary induction on $z$. The rest follows by completeness (Theorem \ref{completeness}) and soundness  (Theorem \ref{soundness}).  
\end{proof}

\begin{lemma}\label{lemma}
	Let $\alpha(z,v)$ be  a progression formula based on $\rho(u,z,v)$. Then 
	\begin{align*} 
		\ha &\vdash  \ \pr_\ha(x)\imp \pr(z,x),  \\
		\ha &\vdash \ \pr(z,x)\imp \pr(2^z,x), \\
		\ha &\vdash  \ \{e\}(n)\simeq z\imp \alpha(z,v)\imp \alpha(3\timess 5^e,v), \\
		\ha & \vdash \ \{e\}(n)\simeq z\imp \pr(z,x)\imp \pr(3\timess 5^e,x).
	\end{align*}
\end{lemma}

\begin{definition}[progression]
Let $\alpha(z,v)$   be a   progression formula based on $\rho(u,z,v)$. 	For $d\in\N$, let 
$T_{d}=\{\psi \mid \alpha(\bar d,\bar \psi) \text{ is true}\}$. We call $(T_d)_{d\in\N}$ a  progression over $\ha$ based on $\rho(u,z,v)$. 	
\end{definition}

\begin{theorem}\label{theorem}
Let $(T_d)_{d\in\N}$ be a  progression over $\ha$ based on the  formula $\rho(u,z,v)$. 	 Then $(T_{d})_{d\in\N}$ is a recursively enumerable sequence of theories such that:	
	\begin{itemize}[leftmargin=5mm]
		\item $\ha=T_0\subseteq T_d$ for every $d$;
		\item $a<_{\ko} b$ implies $T_a\subseteq T_b$;
		\item $T_{3\timess 5^e}=\bigcup_{d<_\ko 3\timess 5^e} T_d$ whenever $3\timess 5^e\in\ko$. 
	\end{itemize}	
\end{theorem}
\begin{proof}
The second item is proved by induction on $<_\ko$. 
\end{proof}

\begin{remark}\label{weak existence}
A close inspection reveals that  we could do without Theorem \ref{soundness}. In fact, 
if $\rho(u,z,v)$ is $\Sigma$ in $\ha$, by the diagonal lemma we can obtain a formula $\alpha(z,v)$ such that $\alpha(z,v)$ is $\Sigma$ in $\ha$ and 
	\begin{align*} 
		\ha \vdash&  \ \alpha_0(v)\imp \alpha(z,v), \\
		\ha \vdash & \  \alpha(0,v) \biimp \alpha_0(v), \\
		\ha \vdash & \  \alpha(2^z,v) \biimp \pr_\ha(\bar \alpha(\dot z,\dot v))\lor \rho(\bar\alpha,z,v),   \\
		\ha \vdash & \ \alpha(3\timess 5^e,v) \biimp \alpha_0(v)\lor \exists n\, \exists z\,  \big(\{e\}(n)\simeq z\land  \pr_{\ha}(\bar \alpha(\dot z,\dot v))\big).
	\end{align*}
The formula $\alpha(z,v)$ would still satisfy  Lemma \ref{lemma} and Theorem \ref{theorem}.
\end{remark}

\begin{definition}[succession formulas for reflection]
Let
\[ \rho(\alpha,z,v)\ =_\text{def}\ \ v=\con^\circ_{\alpha}(\dot z)\]
be the succession formula for consistency.

Let 
\[ \rho(\alpha,z,v)\ =_\text{def}\ v= \exists \vp\, (\st(\vp) \land v= \pr^{\circ}_{\alpha}(\dot z, \dot \vp)\imp^\circ \vp).\]
be the  succession formula for local reflection.

Let
	\[ \rho(\alpha,z,v)\ =_\text{def}\  \exists \vp\, \exists x\, \big(\st(\forall^\circ x\, \vp) \land v=\forall^\circ x\, (\pr^\circ_{\alpha}(\dot z, \dot \vp(\ddot x))\imp^\circ  \vp(x))\big). \]
be the succession formula for uniform reflection.
\end{definition}

\begin{remark}
Let $(T_d)_{d\in\N}$ be a  progression based on the  succession formula for consistency.  Then  $T_{2^d}=T_d+\con(T_d)$ for every $d$. Notice that  $d$ need not be in Kleene's $\ko$. Similarly for local and uniform reflection. 
\end{remark}

\begin{definition}
We say that a formula $\vp$ is provable in an iteration of consistency over $\ha$ if there exists a $d\in\ko$ such that $T_d\vdash \vp$, for some fixed progression $(T_d)_{d\in\N}$ based on the succession formula for consistency.  Similarly for local and uniform reflection. 
\end{definition}

\section{Iterating consistency or local reflection}\label{local sec}

\begin{theorem}[local reflection]\label{local}
Let $\vp$ be a formula. The following are equivalent:
\begin{enumerate}[label=$(\arabic*)$,leftmargin=*]
	\item $\vp$ is provable in an iteration of consistency over $\ha$;
	\item $\vp$ is provable in an iteration of local reflection  over $\ha$;
	\item $\vp$ is provable in $\ha$ $+$ all true $\Pi_1$  sentences.
\end{enumerate}
\end{theorem}
\begin{proof}
Clearly (1) $\Imp$ (2). Let $T$ be the extension of $\ha$  with all true $\Pi_1$ sentences.

(2) $\Imp$ (3) is proved exactly as  in \cite[Thm.\  4.5 p.\ 289]{F62}.  One proves that $T_d\subseteq \{\vp\mid T\vdash \vp\}$ by induction on $d\in\ko$, where $(T_d)_{d\in\N}$ is a progression over $\ha$ based on local reflection. The cases $d=0$ and $d=3\timess 5^e$ are immediate. For the successor case it suffices to show that $T\vdash \pr(\bar d,\bar \vp)\imp \vp$ for every sentence $\vp$  assuming that $T$ proves all theorems of $T_d$. Denote $\pr(\bar d,\bar \vp)$ by $\phi$. The sentence $\neg \phi$ is $\Pi_1$ in $\ha$. We now have two cases. Case (i) $T_d\vdash \vp$. Then by induction, $T\vdash \vp$ and hence $T\vdash \phi\imp \vp$. Case (ii) $T_d\nvdash \vp$. Then $\neg \phi$ is true and hence $T\vdash \neg \phi$. Also in this case we have $T\vdash \phi\imp \vp$. 
Notice that:
\begin{itemize}[leftmargin=5mm]
	\item in case (ii) we use explosion (ex falso) and hence such proof does not extend to, for example, minimal logic;
	\item the case distinction between (i) and (ii), although classical, takes place at the meta level.
\end{itemize}  

(3) $\Imp$ (1) is a verbatim copy of  \cite[Thm.\  4.1 p.\ 287]{F62} once we have Theorem \ref{second}. We reproduce the argument for the sake of the reader. Let $(T_d)_{d\in\N}$ be a progression over $\ha$ based on consistency.  Let $\forall x\ \vartheta(x)$ be a $\Pi_1$ sentence where $\vartheta(x)$ is $\Delta_0$.  By the recursion theorem we can find $e$ such that 
\[  \{e\}(n)\simeq \begin{cases} n_\ko & \text{ if } \forall x\leq n\, \vartheta(x);\\
	2^{3\timess 5^e}  & \text{ otherwise.} \end{cases} \] 
We observe that if  $\forall x\ \vartheta(x)$ is true, then  $d=3\timess 5^e\in\ko$ and so $2^d\in\ko$.  We aim to prove  that $T_{2^d}\vdash \forall x\, \vartheta(x)$, whenever $\forall x\, \vartheta(x)$ is true. Actually, we will show that 
$T_{2^d}\vdash \forall x\, \vartheta(x)$ in any case. This means that if  $\forall x\ \vartheta(x)$ is false, then $2^d\notin\ko$ and $T_{2^d}$ proves a false $\Pi_1$ statement. 

\begin{claim} If  $\exists x\, \neg \vartheta(x)$ is true, then $T_d\vdash \con(\bar d)$. \end{claim}
\begin{proof}[Proof of claim]
	If $\exists x\ \neg \vartheta(x)$ is true,   then $T_{2^d}=T_d$. By construction, $T_{2^d}\vdash \con(\bar d)$. 	
\end{proof}
The claim is provable in $\ha$.
\begin{claim} $\ha\vdash \exists x\, \neg \vartheta(x) \imp \pr(\bar d,\con^\circ(\bar d))$. \end{claim}
\begin{proof}[Proof of claim]
	By the recursion theorem in $\ha$, the definition of $e$ carries over to $\ha$. In particular, 
	\begin{align*} 
		\tag{0} \ \ha &\vdash \{\bar e\}(n)\simeq  2^{\bar d} \biimp \exists x\leq n\, \neg\vartheta(x).
	\end{align*}
	Hence: 
	\begin{align*}	
		\tag{1} \ha &\vdash \{\bar e\}(n)\simeq z \imp \pr(z,x)\imp \pr(3\timess 5^{\bar e},x)	& \text{ Lemma \ref{lemma}}, \\	
		\tag{2} \ha &\vdash \exists x\, \neg \vartheta(x) \imp \pr( 2^{\bar d},x)\imp \pr(\bar d,x)	& \text{ (0) and (1)}, \\
		\tag{3} \ha &\vdash \pr( 2^{\bar d},\con^\circ(\bar d))	& \text{ construction,} \\
		\tag{4} \ha&\vdash 	 \exists x\, \neg \vartheta(x) \imp \pr( \bar d,\con^\circ(\bar d)) & \text{ (2) and (3).}
	\end{align*}
\end{proof}
By Lemma \ref{lemma}, $\alpha(z,v)$ satisfies the requirements of Theorem \ref{second}. This means that 
$\ha\vdash \con(z)\imp \neg\pr(z, \con^\circ(\dot z))$. Thus $\ha\vdash \con(\bar d)\imp \neg\pr(\bar d, \con^\circ(\bar d))$. Since $\ha\subseteq T_{2^d}$ and  $T_{2^d}\vdash \con(\bar d)$, we then obtain $T_{2^d}\vdash \neg\pr(\bar d, \con^\circ(\bar d))$, and so by (4) it follows that $T_{2^d}\vdash \forall x\, \vartheta(x)$, as desired.
\end{proof}

\section{Iterating uniform reflection}\label{uniform sec}

\begin{definition}[recursive $\omega$-rule]
	Let $p_n$ be the $n$-th prime number. Let $\W$ be inductively defined by the following clauses:
	\begin{itemize}[leftmargin=5mm]
		\item if $\vp$ is a logical axiom or an axiom of $\ha$, then $2^{\vp}\in \W$;
		\item if $e$ is the index of a total recursive function, $\{e\}(n)\in \W$ and $\{e\}(n)\simeq 2^{\vp(\bar n)}\timess a_n$ for every $n$, then $2^{\forall x\, \vp(x)} \timess 3 \timess 5^e\in \W$;
		\item if $\vp$ follows from $\vp_0,\ldots,\vp_k$ by means of a logical rule and $a_0,\ldots,a_k\in\W$ with $a_i=2^{\vp_i}\timess b_i$ for all $i\leq k$, then $2^{\vp}\timess 3^2\timess  5^{a}\in\W$, where $a= p_0^{a_0}\cdots p_k^{a_k}$.\footnote{One can choose among the many sets of axioms and rules available in the literature. For example, it is enough to have just two rules, {\em modus ponens} and {\em generalization}, where $k\leq 2$; cf.\ \cite[Ch.\ 2, Sec.\ 4] {TD88}. }
	\end{itemize}
We say that $\vp$ is provable in $\W$ if $2^{\vp}\timess b\in \W$ for some $b$.
\end{definition}

\begin{theorem}[uniform reflection]\label{uniform}
Let $\vp$ be a formula. The following are equivalent:
\begin{enumerate}[label=$(\arabic*)$,leftmargin=*]	
	\item  $\vp$  is provable in an iteration of uniform reflection over $\ha$;
	\item  $\vp$ has   a proof  in $\W$.
\end{enumerate}	
\end{theorem}

We start with the easy implication. From now on, we fix a progression $(T_d)_{d\in\N}$ over $\ha$ based  on uniform reflection. 
\begin{proof}[Proof of $(1) \Imp (2)$]
We claim that there is an index $g$ of a partial recursive function  such that if $d$ is in Kleene's $\ko$ and   $T_d\vdash \vp$,  then $\{g\}(d,\vp)$ is defined and $\{g\}(d,\vp)$ is a proof of $\vp$ in $\W$. The index $g$  can be found with the aid of the recursion theorem. We only indicate how to  look effectively for a proof of 
\[ \forall x\, (\pr(\bar d,\bar \vp(\dot x))\imp\vp(x)), \] 
uniformly in $d$ and $\vp(x)$. 
The search will succeed if $d\in\ko$. Let $\gamma(x)$ be $\pr(\bar d,\bar \vp(\dot x))\imp\vp(x)$. One can find a $\Delta_0$ formula $\vartheta(x,y)$ such that 
\[\ha \vdash \exists y\, \vartheta(x,y)\biimp \pr(\bar d,\bar \vp(\dot x)). \] 
One can therefore construct, effectively in $n$, a proof in $\ha$ and so a proof $c_n$ in $\W$ of 
\[ \forall  y\, (\vartheta(\bar n,y)\imp \vp(\bar n)) \imp \gamma(\bar n). \]
Effectively in $m$, one can decide whether $\vartheta(\bar n, \bar m)$ is true.  If it is true, then $T_d\vdash \vp(\bar n)$. In this case one computes $\{g\}(d,\vp(\bar n))$.  If it is not true, then there is a proof in $\ha$ and so a proof $d_{nm}$ in  $\W$ of $\vartheta(\bar n,\bar m)\imp \vp(\bar n)$. 
So let \[ \{g\}(d,\forall x\, \gamma(x))\simeq 2^{\forall x\, \gamma(x)} \timess 3 \timess 5^e,\] 
 where 
\[ \{e\}(n)\simeq 2^{\gamma(\bar n)}\timess 3^2\timess 5^{a_n}, \ \ \  a_n = 2^{b_n} \timess 3^{c_n}, \ \ \  b_n= 2^{\forall  y\, (\vartheta(\bar n,y)\imp \vp(\bar n))}\timess 3\timess 5^{h_n} \]
and finally
\[ \{h_n\}(m)\simeq \begin{cases}   H_{nm}(\{g\}(d,\vp(\bar n)))   & \text{ if  } \vartheta(\bar n,\bar m); \\
	                      d_{nm}   & \text{ otherwise.} 
	\end{cases} \] 
Here, $H_{nm}$ is a total recursive function such that if $a$ is a proof of $\vp(\bar n)$ in $\W$ then $H_{nm}(a)$ is a proof of  $\vartheta(\bar n,\bar m)\imp \vp(\bar n)$ in $\W$. 
\end{proof}

The hard part is  $(2)\Imp(1)$. The natural approach would be to define a partial recursive function $g$ and then prove by induction on $\W$ that   $g(a)$ is defined, $g(a)\in \ko$  and  $T_{g(a)}\vdash \vp$ whenever $a\in\W$ is a proof of $\vp$. This is exactly how we tackled the forward implication $(1)\Imp(2)$. Let's see where this approach breaks down.   The cases $a=2^\vp$ and $a=2^\vp\timess 3^2\timess 5^b$ can be easily dealt with, by letting $g(a)\simeq 0$ in the first case and  \[ g(a)\simeq g(b_0)+_\ko\cdots +_\ko g(b_k)\] in the second case, where $b=p_0^{b_0}\cdots p_k^{b_k}$. Suppose that $a=2^{\forall x\, \vp(x)}\timess 3\timess 5^e$. By induction on $\W$, one could assume that  \[ \{s\}(n)\simeq  2^{g(\{e\}(0))}+_\ko\cdots +_\ko 2^{g(\{e\}(n))}\in\ko\] and $T_{\{s\}(n)}\vdash \vp(\bar n)$, for every $n$. This would imply $d=3\timess 5^{s}\in\ko$ and  $T_{d}\vdash  \vp(\bar n)$, for every $n$. Now,  it would follow by uniform reflection that $T_{2^d}\vdash \forall x\, \vp(x)$  if there were a proof in $T_{2^d}$ of $\forall x\, \pr(\bar d,\bar\vp(\dot x))$. In that  case, one could let $g(a)\simeq 2^d$.  However, the  induction hypothesis does not  grant us that much.

How do we  fix this?  Simply put, we  avoid induction on $\ko$ and prove everything  inside $\ha$. Indeed, we are able to define a total recursive function  $g$ within $\ha$ (cf.\ Definition \ref{function}) such that
\[ \ha\vdash \forall \vp\, \forall b\, (\fm(\vp)\imp  \pr(g(2^\vp\timess b),\vp)). \] 
The idea, in a nutshell, is that every number $a$ can be regarded as (the code of) a locally correct primitive recursive $\omega$-proof with repetitions of the formula $\vp_a$, where $\vp_a=\vp$ when $a$ is of the form $2^\vp\timess b$. A repetition is an inference of the form \lq\lq $\vp$ follows from $\vp$.\rq\rq\ More in detail, given any number $a$, we can build a primitive recursive tree $S$ (cf.\ Definition \ref{unfolding})  such that  $T_{g(\sigma)}$ proves $\vp_\sigma$ for every $\sigma\in S$, where $\vp_\sigma$ is the formula at node $\sigma$.   The tree $S$ encodes a possibly ill-founded $\omega$-proof with repetitions.  We use repetitions to handle, in particular, the following  cases: 
\begin{enumerate}[label=$\arabic*$.,leftmargin=*]
\item the node $\sigma$ is a pair $2^{\forall x\, \vp(x)}\timess b$ which looks like (the code of) a recursive $\omega$-proof of $\forall x\, \vp(x)$; 
\item  the node $\sigma$ is a pair $2^\vp\timess b$ which does not look like (the code of) a recursive $\omega$-proof of $\vp$; 
\item the node $\sigma$  is a quadruple $\pair{\vp(\bar n),e,s,n}$ and 
$\{e\}(n)$ does not converge in $s$ steps; 
\item the node $\sigma$ is a quadruple $\pair{\vp(\bar n),e,s,n}$ and $\{e\}(n)$ does converge in $s$ steps but not to something that  looks like (the code of) a recursive $\omega$-proof of $\vp(\bar n)$. 
\end{enumerate}

The rest of the section is devoted to carry out the following plan:  
\begin{proof}[Proof plan of $(2) \Imp (1)$]
We  define the index $g$ of a total recursive function and a primitive recursive function $a\mapsto \vp_a$  such that:
%[label=(\roman*),labelindent=-3mm,labelwidth=1.5mm,leftmargin = *]
\begin{enumerate}[label=(1.\arabic*),labelindent=0mm,leftmargin = *]
	\item if $a\in \W$ is a proof of $\vp$ then $\{g\}(a)\in\ko$ and $\vp_a=\vp$;\vspace{1mm}
	\item $T_{\{g\}(a)}\vdash \vp_a$ for every $a$ ($a$ need not be in $\W$).
\end{enumerate} 
\vspace{1mm}
To prove (1.2),  we show that $\ha$ proves the formalized version of (1.2): 
\[ \tag{$\ast$} \ha \vdash \forall a\, \exists d\,  (\{\bar g\}(a)\simeq d\land  \pr(d,\vp_a)). \]
We use Lob's theorem in $\ha$ for $\Pi_2$ sentences (Theorem \ref{lob}) to prove ($\ast$). Namely, we show 
\[ \tag{$\ast\ast$} \ha\vdash \pr_\ha(\overline{\forall a\, \pr(\{\bar g\}(a),\vp_a)}) \imp \forall a\, \pr(\{\bar g\}(a),\vp_a). \]  \end{proof}

Although we can define $+_\ko$ in $\ha$, most of its properties can only be verified by induction on $\ko$, and so are not available in $\ha$.  Crucially, we can establish the following by applying  L\"{o}b's theorem. 

\begin{lemma}\label{right}
	$\ha\vdash \pr(b,x)\imp \pr(a+_\ko b,x)$. 
\end{lemma}
\begin{proof}
	Note that by induction on $\ko$ one can easily prove that $T_b\subseteq T_{a+_\ko b}$, whenever $a+_\ko b\in\ko$.  Moreover, this is true of any recursive progression in the sense of Theorem \ref{theorem}. For this lemma, local reflection is indeed sufficient. The sentence we are trying to prove is $\Pi_2$ in $\ha$.  We can then use L\"{o}b's theorem \ref{lob}. Let  $\theta(c,a,b)$ be the $\Sigma$ formula $c=a+_\ko b$. 
	
	We now reason in $\ha$ and proceed by ordinary induction on $b$. The cases $b=0$ and $b=2^c$ are straightforward: they immediately follow from the induction hypothesis and the definition of $+_\ko$. Suppose $b=3\timess 5^e$ and $\pr(b,x)$. Let $c=a+_\ko b$. Recall that $c=3\timess 5^l$, where $\{l\}(n)\simeq a+_\ko 2^{\{e\}(n)}$. We are going to use the L\"{o}b's hypothesis \[ \tag{$\ast$} \pr_\ha(\overline{\forall a\, \forall b\, \forall c\,  \forall x\, (\pr(b,x)\land \theta(c,a,b) \imp \pr(c,x))}). \]
	By definition, $\pr(b,x)$ means that there is a proof  $y$  of $x$ such that  $\alpha(b,y_i)$ or $\infer(y_i,y\restriction  i)$ for every $i<\lh(y)$. We claim that $\pr(c,y_i)$ for every $i<\lh(y)$. If we prove this, then we are done. It is sufficient to show that  $\alpha(b,x)$ implies $\pr(c,x)$. By construction, if $\alpha(b,x)$ then $\{e\}(n)\simeq z$ and $\alpha(z,x)$ for some $n$ and $z$. In particular, $\pr(z,x)$. Let $d=a+_\ko z$.  By $\Sigma$ completeness, $\pr_\ha(\pr^\circ(\dot z,\dot x)\land \bar \theta(\dot d,\dot a,\dot z))$. It follows  by ($\ast$) that $\pr_\ha(\pr^\circ(\dot d,\dot x))$. By Lemma \ref{lemma},  $\pr(2^d,\pr^\circ(\dot d,\dot x))$, and hence by (local) reflection $\pr(2^d, x)$. By definition,  $\{l\}(n)\simeq a+_\ko 2^z= 2^{a+_\ko z}=2^d$. Thus $\pr(\{l\}(n), x)$. By Lemma \ref{lemma}, we finally obtain $\pr(c,x)$, as desired. 
\end{proof}

\begin{definition}
	There is a primitive recursive function $L$ such that 
	\begin{align*} \{L(e)\}(n) &\simeq 2^{\{e\}(0)}+_\ko\cdots +_\ko 2^{\{e\}(n)}.
	\end{align*}
	Let $G(e)=3\timess 5^{L(e)}$. Note that the definition of $G$ carries over to $\ha$.
\end{definition}

\begin{lemma}\label{G}
	Let $G$ be as above. Then the following hold:
	\begin{enumerate}[label=\upshape{(\roman*)},labelindent=-3mm,leftmargin=*]
		\item if $\{e\}(n)\in\ko$ for every $n$, then $G(e)\in \ko$;
		\item $\ha\vdash (\forall n\, \{e\}(n)\downarrow)  \imp  (\{e\}(n)\simeq z \land \pr(z,x))\imp \pr(G(e),x)$.
	\end{enumerate}	
\end{lemma}
\begin{proof}
	(i) is immediate. For (ii), use Lemma \ref{lemma} and  Lemma \ref{right}.  
\end{proof}

\begin{notation}
Let $\pair{\vp,e,m,z}$ denote  $2^{\vp}\timess 3^3\timess 5^{e}\timess 7^{\, 2^m\timess 3^z}$.  	
\end{notation}

\begin{definition}[the unfolding of a proof]\label{unfolding}
Let $S$  be the primitive recursive function defined by:
\[    S(a,n)   \ = 
\begin{cases}
	0          & \text{if } a=2^{\vp} \text{ and } \vp \text{ is a logical axiom or an axiom of } \ha \\&   \text{ ($a$ is axiomatic);}  \\                                           
 \pair{ \vp(\bar n),e,n,0} &  \text{if } a=2^{\forall x\, \vp(x)} \timess 3 \timess 5^e \text{ ($a$ is universal);}  \\ 
      \pair{\vp,e,m,z+1}    &  \text{if } a= \pair{\vp,e,m,z}   \text{ and } \{e\}_z(m) \text{ is undefined;} \\ 
	  b             &  \text{if }  a= \pair{\vp,e,m,z}  \text{ and } \{e\}_z(m)\simeq b \text{ and }  b= 2^\vp\timess c;\\
	 2^\vp                 & \text{if } a= \pair{\vp,e,m,z}  \text{ and } \{e\}_z(m)\simeq b \text{ and }  b \text{ is not} \\ & \text{ of the form } 2^\vp\timess c; \\
		a_n  & \text{if } a=2^{\vp}\timess 3^2\timess  5^{a} \text{ and } a=p_0^{a_0}\cdots p_k^{a_k} \text{ with }  \\ & \ a_i=2^{\vp_i}\timess b_i \text{ for every $i\leq k$ and } \vp \text{ follows from } \\ & \text{ $\vp_0,\ldots,\vp_k$ by means of a logical rule and } n\leq k  \\
	& \text{ ($a$ is an inference of arity $k$);}     \\
	2^\vp        &   \text{if }  a=2^\vp\timess b \text{ and none of the above hold;} \\
	0      & \text{if none of the above holds.}
\end{cases} \]
We say that $a$ is a repetition if $a$ is neither axiomatic nor universal nor an inference. Finally,  
let \[ \vp_a= \begin{cases} \vp & \text{ if } a=2^\vp\timess b; \\
	  \bot & \text{ otherwise.} \end{cases} \]
\end{definition}

\begin{definition}\label{function}
By the recursion theorem we can find $g$ such that 
\[ \{g\}(a) \simeq 
\begin{cases}
0 & a \text{ is axiomatic;} \\
2^{G(H(g,a))} & a   \text{ is universal;} \\
2^{G(I(g,a))}  &   a \text{ is an inference of arity k;} \\
2^{3\timess 5^{J(g,a)}}   & a \text{ is a repetition.}  
\end{cases} \]
Here,
\begin{align*}
 \{H(g,a)\}(n) &\simeq \{g\}(S(a,n)),  \\
 \{I(g,a)\}(n) & \simeq \begin{cases} \{g\}(S(a,n))  &  \text{if }  n\leq k;\\
 	0 & \text{otherwise,}
 \end{cases} \\
\{J(g,a)\}(n)& \simeq \{g\}(S(a,0))+_\ko n_\ko.
\end{align*}
Notice that $g$ is the index of a total recursive function.
\end{definition}

\begin{lemma}\label{i}
Let $g$ be as above. If $a\in\W$, then $\{g\}(a)\in\ko$. 
\end{lemma}
\begin{proof}
The proof is by induction on $\W$. We first observe that if $a\in\W$, then  $a$ is either axiomatic,  universal or an inference. The cases where $a$ is axiomatic or an inference are straightforward. Suppose $a$ is universal, say $a= 2^{\forall x\, \vp(x)} \timess 3 \timess 5^e$. By Lemma \ref{G} part (i) and the definition of $g$, it suffices to show that  $\{g\}(S(a,n))\in\ko$ for every $n$. Fix $n$ and let  $a_{z}= \pair{ \vp(\bar n),e,n,z}$. Since $a\in\W$, we have that  $\{e\}(n)\simeq b$ for some $b\in\W$. By induction on $\W$, we can assume $\{g\}(b)\in\ko$. Now let $\hat z$ be least (indeed the unique $z$) such that $\{e\}_{z}(n)\simeq b$. Then, for any given $m$,  $S(a_{\hat z},m)\simeq b$ and   $S(a_{z},m)=a_{z+1}$ for $z<\hat z$.  By definition, 
\[ \{g\}(a_{\hat z})\simeq 2^{3\timess 5^{J(g,a_{\hat z})}}, \text{ where } \{J(g,a_{\hat z})\}(m) \simeq \{g\}(b)+_\ko m_\ko \]
and for $z<\hat z$ 
\[ \{g\}(a_{z})\simeq  2^{3\timess 5^{J(g,a_{z})}}, \text{ where } \{J(g,a_{z})\}(m) \simeq \{g\}(a_{z+1})+_\ko m_\ko. \]
By an ordinary induction (on $\hat z-z$),\footnote{We find echoes of such  \lq\lq backwards induction\rq\rq\ in  Feferman, where he describes his proof idea in terms of \lq\lq backwards recursion\rq\rq; cf.\ \cite[p.\ 293]{F62}.} we can then prove that $\{g\}(a_{z})\in\ko$ for every $z\leq \hat z$. But $S(a,n)=a_{0}$. Hence  $\{g\}(S(a,n))\in\ko$. This concludes the proof. 
\end{proof}

Finally, we can deliver on our promise. 

\begin{proof}[Proof of $(2)\Imp(1)$]
By Lemma \ref{i}, which proves (1.1), we need only show \vspace{1mm}

\begin{enumerate}[label=(1.\arabic*),labelindent=0mm,leftmargin = *,start=2] 
\item $T_{\{g\}(a)}\vdash \vp_a$  for every $a$.
\end{enumerate} \vspace{1mm}

The index $g$ can be defined by the recursion theorem and proved total within  $\ha$.  Drawing from our earlier discussion (outlined in the proof plan), it suffices to prove in $\ha$  the following claim. 
\vspace{1mm}
\begin{enumerate}[label=(1.\arabic*),labelindent=0mm,leftmargin = *,start=3]
\item if $\ha\vdash \forall a\, \pr(\{\bar g\}(a),\vp_a)$, then $T_{\{g\}(a)}\vdash \vp_a$ for every $a$.
\end{enumerate}
\vspace{1mm}
The details of the formalization are left to the reader. 
\begin{proof}[Proof of claim $(1.3)$]
Assume $\ha\vdash \forall a\, \pr(\{\bar g\}(a),\vp_a)$ and let $a$ be given. We proceed by cases. \vspace{1mm}

Case 1. $a=2^\vp$ is axiomatic. Then  $\vp_a=\vp$ is a logical axiom or an axiom of $\ha$   and $\{g\}(a)\simeq 0$. We can directly conclude that $T_{\{g\}(a)}\vdash \vp_a$ since  $T_{\{g\}(a)}=T_0=\ha$ and $\ha\vdash \vp$. \vspace{1mm}

Case 2. $a$ is universal. Say $a=2^{\forall x\, \vp(x)} \timess 3 \timess 5^e$. In this case, $\vp_a=\forall x\, \vp(x)$ and  $\vp_{S(a,n)}=\vp(\bar n)$ for all $n$.  Moreover, $\{g\}(a)\simeq 2^{G(H(g,a))}$, where  $\{H(g,a)\}(n) \simeq \{g\}(S(a,n))$.  All this formalizes in $\ha$. In particular,  
\begin{align*}
\ha &\vdash \forall x\, (\vp_{S(\bar a,x)}=\bar \vp(\dot x)), \\
\ha &\vdash  \forall x\, \exists z\, (\{\bar g\}(S(\bar a,x))\simeq z\land  \pr(z,\vp_{S(\bar a,x)})), \\
\ha & \vdash \forall x\, (\{H(\bar g,\bar a)\}(x) \simeq \{\bar g\}(S(\bar a,x))),
\end{align*}
where the second item follows from the assumption. Then by Lemma \ref{G} part (ii)
\[ \ha\vdash \forall x\, \pr(G(H(\bar g,\bar a)),\bar \vp(\dot x)). \]
Let $d=G(H(g,a))$ so that $\{g\}(a)\simeq 2^d$. Since  $\ha\subseteq T_{2^d}$ and $\ha\vdash G(H(\bar g,\bar a))=\bar d$, we obtain
\[ \tag{1} T_{2^d}\vdash \forall x\, \pr(\bar d,\bar \vp(\dot x)). \]
On the other hand, by definition,
\[  T_{2^d}\vdash \forall x\, (\pr(\bar d,\bar \vp(\dot x)) \imp \vp(x)). \]
Hence $T_{2^d}\vdash \forall x\, \vp(x)$, as desired. If there are  free variable in $\forall x\, \vp(x)$, say $z_1,\ldots,z_m$, take the universal closure $\psi(x)=_{def} \forall z_1\cdots \forall z_m\, \vp(x,z_1,\ldots,z_m)$. Now, from (1) we can infer $T_{2^d}\vdash \forall x\, \pr(\bar d,\bar \psi(\dot x))$, and then conclude $T_{2^d}\vdash \forall x\, \psi(x)$ by applying uniform reflection with respect to $\psi(x)$. In particular, $T_{2^d}\vdash \forall x\, \vp(x)$. 
    \vspace{1mm}

Case 3. $a$ is an inference of arity $k$. Say $a=2^{\vp}\timess 3^2\timess  5^{a}$ and $a= p_0^{a_0}\cdots p_k^{a_k}$  with $a_i=2^{\vp_i}\timess b_i$ for $i\leq k$. In this case, $\vp_a=\vp$ and $\vp_{S(a,n)}=\vp_n$ for every $n\leq k$. Moreover, $\{g\}(a)\simeq 2^{G(I(g,a))}$, where $\{I(g,a)\}(n) \simeq \{g\}(S(a,n))$ for all   $n\leq k$.  As before, everything formalizes in $\ha$. In particular, for all $n\leq k$,
\begin{align*}
	\ha &\vdash \vp_{S(\bar a,\bar n)}=\overline{\vp_n}, \\
	\ha &\vdash  \exists z\, (\{\bar g\}(S(\bar a,\bar n))\simeq z\land  \pr(z,\vp_{S(\bar a,\bar n)})), \\ 
		\ha &\vdash  \{I(\bar g,\bar a)\}(\bar n) \simeq \{\bar g\}(S(\bar a,\bar n)).
\end{align*}
Then by Lemma \ref{G} part (ii),   for all $n\leq k$,
\[ \ha \vdash  \pr(G(I(\bar g,\bar a)),\overline{\vp_n}). \]
Let $d=G(I(g,a))$ so that $\{g\}(a)\simeq 2^d$. As before, we obtain 
\[ T_{2^d}\vdash \pr(\bar d,\overline{\vp_n}) \]
for all $n\leq k$. Since uniform reflection implies local reflection, we thus have 
\[ T_{2^d}\vdash \vp_n \]
for all $n\leq k$. Finally, since $\vp$ follows from $\vp_0,\ldots,\vp_k$ by means of a logical rule, we get $T_{2^d}\vdash \vp$, as desired. \vspace{1mm}

Case 4. $a$ is a repetition. Then $\vp_a=\vp_{S(a,0)}=\vp$.  This includes the case $\vp=\bot$. Moreover, $\{g\}(a)\simeq 2^{3\timess 5^{J(g,a)}}$, where $\{J(g,a)\}(n)\simeq \{g\}(S(a,0))+_\ko n_\ko$. Then
\begin{align*}
		\ha &\vdash \vp_{S(\bar a,0)}=\bar{\vp}, \\
		\ha &\vdash \exists z\, (\{\bar g\}(S(\bar a,0))\simeq z\land  \pr(z,\vp_{S(\bar a,0)})), \\ 
		\ha &\vdash  \{J(\bar g,\bar a)\}(\bar n) \simeq \{\bar g\}(S(\bar a,0))+_\ko n_\ko. 
\end{align*}
In particular,  \[ \ha\vdash \pr(\{J(\bar g,\bar a)\}(0),\bar \vp) \]
and therefore \[ \ha\vdash \pr(3\timess 5^{J(\bar g,\bar a)},\bar \vp). \] 
Let $d=3\timess 5^{J(g,a)}$ so that $\{g\}(a)\simeq 2^d$. Then 
\[ T_{2^d}\vdash \pr(\bar d, \bar \vp). \] 
Again, since uniform reflection implies local reflection, we have $T_{2^d}\vdash \vp$, as desired. \end{proof}
%Our attentive reader can check that the  proof of the claim can be carried out in $\ha$.
\end{proof}

%\bibliographystyle{plain}
%\bibliography{biblio}{}

\end{document}